*We write here, without a proof, a beautiful and remarkable proposition that we discovered: ... I don't believe that it can be given on numbers a more beautiful and general theorem. But I have no time and no space to bring back a proof on the margin.* Pierre De Fermat [1]

# A GEOMETRICAL WAY TO SUM POWERS BY MEANS OF TETRAHEDRONS AND EULERIAN NUMBERS

MARIO BARRA

ABSTRACT. We geometrically prove that in a d-dimensional cube with edges of length n, denoted with $C^d{}_n$, we have $C^d{}_n = n^d = \sum_{s=0}^{d-1} \left\langle {d \atop s} \right\rangle T^d_{n-s}$,[2] $\forall n,d,s \in N$, where the Eulerian number[3] $\left\langle {d \atop s} \right\rangle$ becomes the number of particular d-dimensional tetrahedrons with edges of length n-s,[4] $T^d{}_{n-s}$, that tessellate $C^d{}_n$. As a consequence, in order to add powers of integers, 1, 2 or t times, applying the well known property $T^d{}_n = \sum_1^n T^{d-1}{}_h$, 1, 2 or t times, we have

$$\sum_{j=1}^n j^d = \sum_{s=0}^{d-1} \left\langle {d \atop s} \right\rangle T^{d+1}_{n-s}, \quad \sum_{m=1}^n \sum_{j=1}^m j^d = {}^2\sum_{j=1}^n j^d = \sum_{s=0}^{d-1} \left\langle {d \atop s} \right\rangle T^{d+2}_{n-s}, \quad {}^t\sum_{j=1}^n j^d = \sum_{s=0}^{d-1} \left\langle {d \atop s} \right\rangle T^{d+t}_{n-s}.$$

Finally we prove that $\left\langle {d \atop s} \right\rangle = \nabla^{d+1}(s+1)^d$.

## 1. INTRODUCTION AND EXAMPLES

We want to calculate the numbers of tetrahedrons with equal dimension and decreasing edges, that operate a dissection of a d-dimensional cube with edges of length n, $C^d{}_n$. Such tetrahedrons are built with the solutions of systems with the signs " $\geq$" and " > " between the variables $X_1,X_2,\ldots,X_d$. Examples, for d=3, $X_1 \geq X_2 \geq X_3$ e $X_3 > X_1 \geq X_2$.

Maybe it will be better to introduce the argument through some examples.

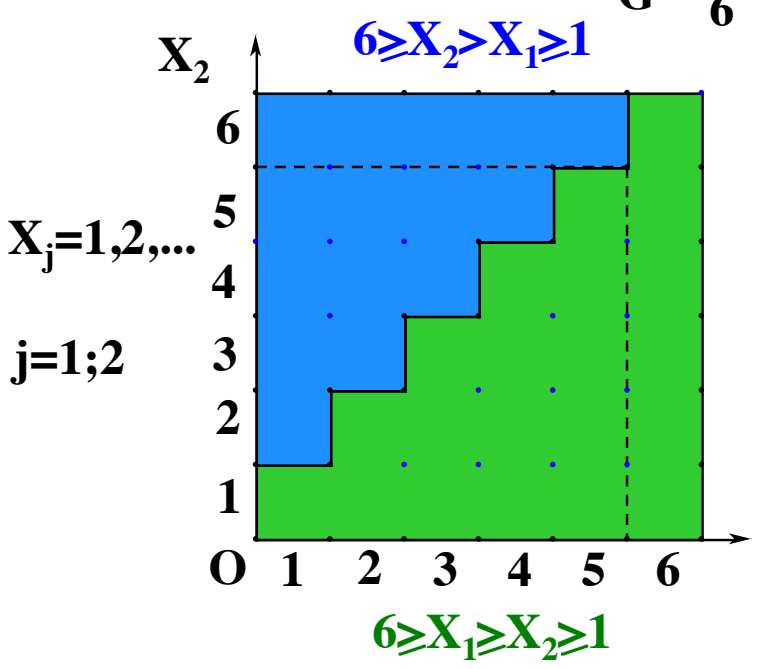

$$G^{2\wedge}{}_6 = C^2{}_6 - C^2{}_5 \qquad C^2{}_6 = \sum_{h=1}^{6} G^{2\wedge}{}_h$$

$$G^{2\wedge}{}_6 = T^1{}_6 + T^1{}_5$$

$$C^2{}_6 = T^2{}_6 + T^2{}_5$$

$$\sum^6 C^2{}_h = 36 + 25 + 16 + 9 + 4 + 1 = 91 = T^3{}_6 + T^3{}_5 = 56 + 35$$

$$\sum_{m=1}^{6} \sum^m C^2{}_h = 1 + 5 + 14 + 30 + 55 + 91 = 196 = T^4{}_6 + T_5{}^4 = 126 + 70$$

...

**Example 1**.

---

[1] Pierre De Fermat, *Osservazioni su Diofanto, Boringhieri*, 1959 (1670), p. 108. Fermat speaks about the formula that expresses the number of points of a tetrahedron in d dimensions with edges of length n, equivalent to: $T^d{}_n = \frac{n(n+1)\cdots(n+d-1)}{d!}$. In the same book, p. 18, speaking about its *Last Theorem*, Fermat writes in a more measured way: ... *it is impossible to divide a cube into two cubes, ..., or generally any power beyond the square into like powers; of this I have found a remarkable demonstration. This margin is too narrow to contain it.*

[2] The sign "=" after the symbol of a polyhedron indicates its cartesian definition, or its measure, or the formula to calculate the measure. Other meanings of symbols: $n^{j\geq} = n(n+1)\ldots(n+j-1)$, $n^{j>} = (n+1)(n+2)\ldots(n+j)$.

[3] The Eulerian numbers: $\left\langle {d \atop s} \right\rangle$, d≥0, are defined by: $\left\langle {d \atop s} \right\rangle = (1+s)\left\langle {d-1 \atop s} \right\rangle + (d-s)\left\langle {d-1 \atop s-1} \right\rangle$, $\left\langle {1 \atop 0} \right\rangle = 1$, $\left\langle {d \atop k} \right\rangle = 0$, k≥d.

[4] The word "edge" also indicates the side of a triangle. In general, the terms used in the third dimension will be used in any dimension.

$C^2_6 = 6^2 = \mathbf{1}T^2_6 + \mathbf{1}T^2_5 = \left\langle{2 \atop 0}\right\rangle T^2_6 + \left\langle{2 \atop 1}\right\rangle T^2_5$, where $\left\langle{2 \atop 0}\right\rangle = \mathbf{1}$, and $\left\langle{2 \atop 1}\right\rangle = \mathbf{1}$ indicate that $C^2_6$ is divided respectively in: **1** tetrahedron of dimensions 2, $T^2_6$, with edges of length 6, and **1** tetrahedrons of dimensions 2, $T^2_5$, with edges of length 5.

**Example 2.** For d = 3, trying to imagine what is not visible in the rest of the cube

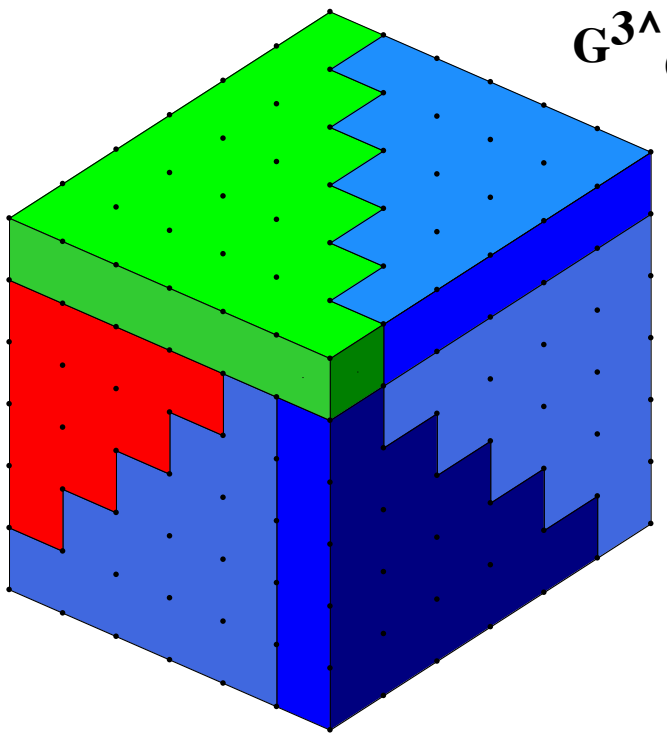

$$G^{3\wedge}_6 = C^3_6 - C^3_5 = 91 \quad C^3_6 = \sum_1^6 G^{3\wedge}_h$$

$$G^{3\wedge}_6 = T^2_6 + 4\,T^2_5 + T^2_4 = 91$$

$$C^3_6 = T^3_6 + 4\,T^3_5 + T^3_4 = 56 + 140 + 20 = 216$$

$$\sum_1^6 C^3_h = T^4_6 + 4\,T^4_5 + T^4_4 = 441 = 126 + 280 + 35$$

$$\sum_{m=1}^{6} \sum_{h=1}^{m} C^3_h = T^5_6 + 4\,T^5_5 + T^5_4 = 812. \quad \dots$$

$C^3_6 = 6^3 = \mathbf{1}T^3_6 + \mathbf{4}T^3_5 + \mathbf{1}T^3_4 = \left\langle{3 \atop 0}\right\rangle T^3_6 + \left\langle{3 \atop 1}\right\rangle T^3_5 + \left\langle{3 \atop 2}\right\rangle T^3_4$, where $\left\langle{3 \atop 0}\right\rangle = \mathbf{1}$, $\left\langle{3 \atop 1}\right\rangle = \mathbf{4}$, $\left\langle{3 \atop 2}\right\rangle = \mathbf{1}$, indicate that $C^3_6$ is divided respectively in: **1** tetrahedron of dimension 3, with edges of length 6, $T^3_6$, **4** tetrahedrons of dimensions 3, with edges of length 5, $T^3_5$, and **1** tetrahedron of dimensions 3, with edges of length, 4, $T^3_4$. The length of the edges are indicated below.

## 2. DEFINITIONS AND SIMPLE PROPERTIES

Some definitions and simple properties.

- **Tetrahedrons and tetrahedral numbers**. We have a tetrahedrons $\mathbf{T^d(n)^+}$ if the *rule of the sections* is valid: $\mathbf{T^d(n)^+} = \bigcup_0^n \mathbf{T^{d-1}(h)^+}$, where $T^0(h)^+$, $\forall \mathbf{h}$, contains only one unit (a point or a little cube) and where $\mathbf{T^{d-1}(h)^+}$, h=1,2,…,n, are disjointed and have h+1 unities on every edges. Therefore we have $\mathbf{T^d_{n+1}} = \sum_1^{n+1} \mathbf{T^{d-1}_h}$, $\forall$d, n$\in\mathbf{N}$, where $\mathbf{T^d_{h+1}}$ is the measure of all the types of tetrahedrons containing h+1 unities on every edge.

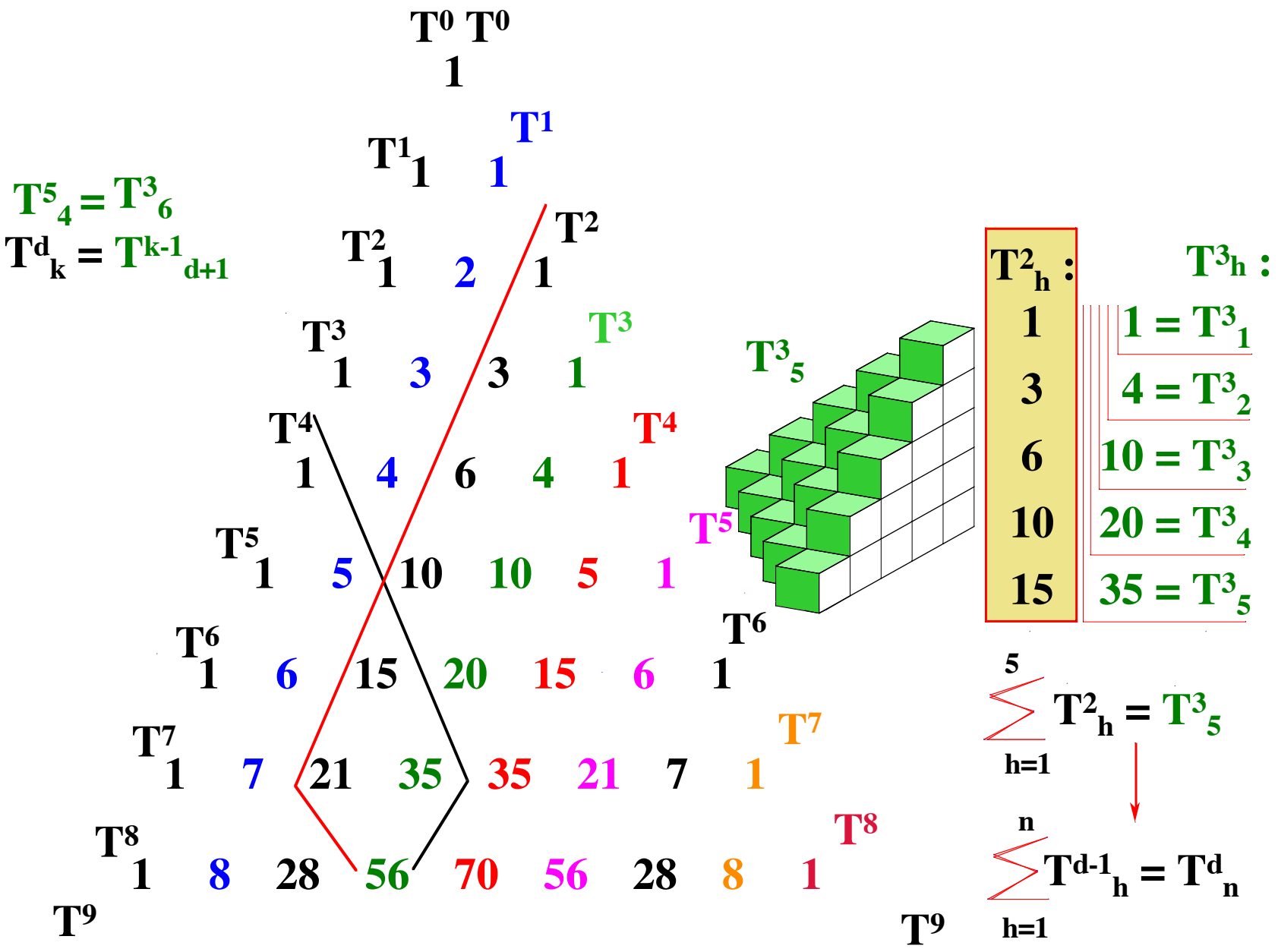

**- Zone-coordinates**: The coordinates $X_h=k$, $k=1,2,...$, $h=1,2,...,d$, are equivalent to $k-1\le x_h\le k$, $x_h\in\mathbf{R}$. $X_h\in\underline{\mathbf{N}}^d$, where $\underline{\mathbf{N}}^d$ is the positive portion of space tessellated with unitary little cubes $c^d=(X_1,X_2,...,X_d)$ whose edges are parallel to the coordinate axes and whose faces are in common with the contiguous [5] little cubes.

**- Fishbone**: Particular system of equations or inequalities constituted with a sequence of variables separated by the signs of inequality or equality that indicate which zone-coordinates verify the system.

- **Euler's fishbone**: Fishbone with only the signs "≥" and ">" between the variables, established with the **rule of the climbs** that places between every couple of variables $X_h$, $X_k$ the sign "≥" only if k>h, and, otherwise, it inserts the sign ">". Examples: $X_1\ge X_2\ge X_3$ and $X_3>X_1\ge X_2$.

**- Extremes of a fishbone**: Integers, m, M, m≤M, that indicate the maximum and the minimum of the values that the variables can assume. Example: $\mathbf{4}\ge X_1\ge X_2\ge X_3\ge\mathbf{1}$, that is equivalent to consider $X_1\ge X_2\ge X_3$ in $C^3{}_4=\{(X_1,X_2,X_3): X_h=\mathbf{1},2,3,4,\ h=1,2,3\}$.

**- Euler's Escher**: Set of *Euler's fishbones* constructed taking all the d! permutations $X_{h_1}, X_{h_2},..., X_{h_d}$ of the variable $X_1, X_2,..., X_d$, and placing between these the signs " > " and " ≥" according to the *rule of the climbs*.

- **Two properties of an Escher**: Given an Escher and the coordinates of any little cube $c^d$ in $C^d{}_n$, we have the following properties:

**1)** It is possible to establish the *fishbone* satisfied by the coordinates of $c^d$ observing only the order of the coordinates.

In fact, if we order in a non-increasing order the numerical coordinates of $c^d$, putting $X_h$ before $X_k$, if $X_h=X_k$ and h<k, these ones satisfy only the *fishbone* expressed with the correspondent variables in the same order, placing between these the signs ">" and "≥" according to the *rule of the climbs*.

**2)** The coordinates of any $c^d$ can satisfy only one *fishbone*.

In fact the specified order is unique.

Therefore an Escher operates a partition of $C^d{}_n$.

A different proof of the last achievement will be introduced after the proof that $\mathbf{n^d = \sum_{s=0}^{d-1}\left\langle {d \atop s} \right\rangle T^d_{n-s}}$.

**Example 3**.

The coordinates of any $c^3$ in $\mathbf{C^3{}_4}$, $c^3=(3,2,4)$, watching the *rule of the climbs*, can only satisfy the *fishbone* $X_3>X_1\ge X_2$. In a similar way (2,2,4) satisfy only the previous *fishbone* $X_3>X_1\ge X_2$, and $c^3=(6,6,1,6,3,1,5,9)$ satisfy only the *fishbone* $X_8>X_1\ge X_2\ge X_4\ge X_7>X_5>X_3\ge X_6$.

Once the cubes have been tessellated into tetrahedrons, the sum of the powers of integers becomes banal: it is enough to apply the property $\mathbf{T^d{}_n=\sum_1^n T^{d-1}{}_h}$.

[5] Two little cubes are contiguous when their coordinates have unitary Hamming's distance (all the coordinates are equal, except one of them that differs by an unit).

## 3. REGULAR AND RECTANGULAR TETRAHEDRONS

**Lemma 1.** *The following sets, A and B, represent two tetrahedrons, called respectively regular tetrahedron,* $T^d(n)^+$, *and rectangular tetrahedron,* $T^d(\leq n)^+$.

$A=\{(x_1,x_2,...,x_{d+1})\colon \sum_{1}^{d+1} \mathbf{x_j} = n,\ x_j= \mathbf{0},1,2, ...,\ j=1,2,...,d+1\}$,

$B=\{(x_1,x_2,...,x_d)\colon \sum_{1}^{d} \mathbf{x_j} \leq n,\ x_j=\mathbf{0},1,2,...,\ j=1,2,...,d\}$.

*Proof.* The sets A and B are tetrahedrons because they verify the *rule of the sections*. In fact $A = \bigcup_{h=0}^{n} \{(x_1,x_2,...,x_{d+1})\colon \sum_{j=1}^{d} \mathbf{x_j} = n-h,\ x_j=\mathbf{0},1,2, ...,\ j=1,2,...,d,\ x_{d+1}=h\}$.

$B = \bigcup_{h=0}^{n} \{(x_1,x_2,...,x_d)\colon \sum_{j=1}^{d-1} \mathbf{x_j} \leq n-h,\ x_j=\mathbf{0},1,2,...,\ j=1,2,...,d-1,\ x_d=h\}$.

**Lemma 2.** *The regular tetrahedron*: $T^d(n)^+=\{(x_1,x_2,...,x_{d+1})\colon \sum_{1}^{d+1} \mathbf{x_j}=n,\ x_j= \mathbf{0},1,2, ...,$ $j=1,2,...,d+1\}$ *has the following number of points*: $T^d(n)^+= \frac{\mathbf{n^{d>}}}{\mathbf{d!}} = \binom{\mathbf{n+d}}{\mathbf{d}} =$ $= \frac{\mathbf{(n+1)\cdots(n+d)}}{\mathbf{d!}} = T^d_{n+1} = T^n(d)^+$. $T^d(n)^+$ *has n+1 points on every edge.*

*Proof.* It is sufficient to generalize what follows. The number of all the anagrams of the word composed with 3 dots and 2 bars: " ... | | " is $\frac{5!}{3!2!} = \binom{3+2}{2} = \binom{5}{2} = 10$. They are

| ...\| \| | ..\|.\| | ..\| \|. | .\|..\| | .\|.\|. | .\| \|.. | \|...\| | \|..\|. | \|.\|.. | \| \|... |
|---|---|---|---|---|---|---|---|---|---|
| 3+0+0 | 2+1+0 | 2+0+1 | 1+2+0 | 1+1+1 | 1+0+2 | 0+3+0 | 0+2+1 | 0+1+2 | 0+0+3 |

Interpreting the number of dots:

- before the first bar, as the value of $x_1$,
- between the first and second bar, as the value of $x_2$,
- after the second bar, as the value of $x_3$,

all the ways are obtained in order to have sum 3 with 3 addends. Then

$\mathbf{T^2(3)^+}=\{(x_1,x_2,x_3)\colon \sum_{1}^{3} \mathbf{x_j}=3,\ x_j=\mathbf{0},1,2,...,\ j=1,2,3\}=\frac{\mathbf{3^{2>}}}{\mathbf{2!}}=\binom{\mathbf{3+2}}{\mathbf{2}}=\frac{\mathbf{4\cdot 5}}{\mathbf{2!}}=T^2_4=\binom{\mathbf{3+2}}{\mathbf{3}}=\mathbf{T^3(2)^+}=T^3_3$

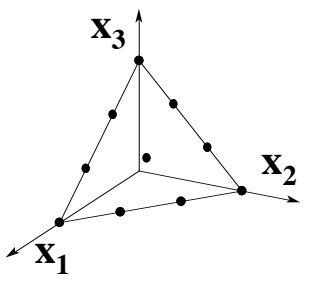

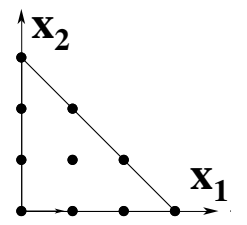

$\mathbf{T^2_4}=\frac{\mathbf{3^{2>}}}{\mathbf{2!}}=\frac{\mathbf{4^{2\geq}}}{\mathbf{2!}}=\frac{\mathbf{4\cdot 5}}{\mathbf{2!}}=\frac{\mathbf{3^{2>}}}{\mathbf{2!}}=\mathbf{T^3_3}$

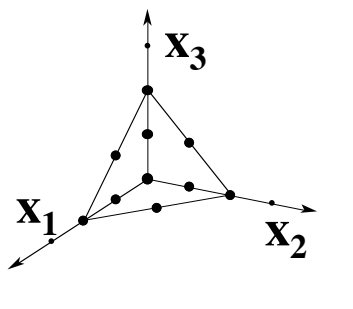

**Lemma 3.** *The regular tetrahedron*, $T^d(n)^+$, *has the same points of the rectangular tetrahedron,* $T^d(\leq n)^+$.

*Proof.* It is sufficient to generalize what follows. In order to determine the points of $T^2(\leq 3)^+=\{(x_1,x_2)\colon \sum_{1}^{2} \mathbf{x_j} \leq 3,\ x_j= \mathbf{0,1,2,} ...,\ j=1,2\}$, it can be considered that:

1) to have sum 3 with 3 addends, 3 dots and 2 bars are necessary;

2) to have sum 3 with 2 addends, 3 dots and 1 bars are necessary. We add in this second case at the end one other bar for having after, in the anagrams, from zero to 3 dots, so as to have before a number $\leq 3$ of dots and therefore a sum $\leq 3$, in all the ways. In both cases 3 dots and 2 bars are necessary. In general terms, as for sum n with d+1 addends, as for sum $\leq n$ with d addends, the dots and the bars we need are always n and d: the anagrams of the word with these "letters" prove that:

$$T^d(n)^+=T^d(\leq n)^+=\binom{\mathbf{n+d}}{\mathbf{d}}=\frac{\mathbf{(n+1)\cdots(n+d)}}{\mathbf{d!}}=\frac{\mathbf{n^{d>}}}{\mathbf{d!}}=\frac{\mathbf{(n+1)^{d\geq}}}{\mathbf{d!}}= T^d_{n+1}= \binom{\mathbf{d+n}}{\mathbf{d}}= T^n(d)^+ = T^n_{d+1}$$

**Example 4.**

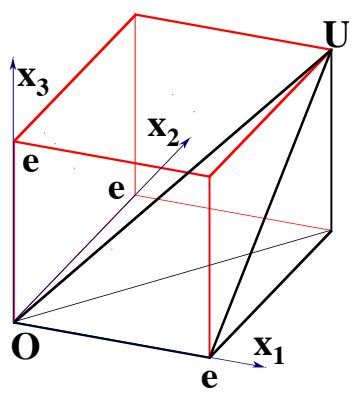

$e \geq x_1 \geq x_2 \geq x_3 \geq 0$

- The cube $C^3(e) = \{(x_1,x_2,x_3)$: $0 \leq x_h \leq e$, $h=1,2,3\} = e^3$, $e \in \mathbf{R}$, is divided in 3! tetrahedrons: $T!^3(e) = \{(x_i,x_j,x_h)$: $e \geq x_i \geq x_j \geq x_h \geq 0$, ... $i \neq j \neq h\}$, where i,j,h is one of the 3! permutations of indices 1,2,3. All the tetrahedrons, for symmetry, have the same measure. Since $e \geq x_h \geq 0$, h=1,2,3, $T!^3(e)$ has 3 edges of length **e**, that form a broken line, respectively parallels to coordinate axes: $x_i,x_j,x_h$. It follows: $T!^3(e) = \frac{e^3}{3!}$, and, in general: $T!^d(e) = \frac{e^d}{d!}$.[6]

- The cube: $C^3_n = \{(X_1,X_2,X_3)$: $X_h = 1,2,...,n$, $h=1,2,3\} = n^3$, $n \in \mathbf{N}$, contains the 3! tetrahedrons: $T^3_n = \{(X_i,X_j,X_h)$: $n \geq X_i \geq X_j \geq X_h \geq 1, i \neq j \neq h\} = \frac{n(n+1)(n+2)}{3!} = \frac{n^{3\geq}}{3!} = T^3_n$, where: i,j,h is a permutation of 1,2,3.

Since: $n \geq X_h \geq 1$, $\forall h$, it follows that every $T^3_n$ contains n little cubes on every edge, three of which are respective parallels to coordinate aces $X_i,X_j,X_h$ and form a broken line.

Naturally $C^3_n$ contains also other tetrahedrons. Some of these are characterized by inequalities analogous to the previous one, where some signs " ≥" are replaced with the sign ">". Many of these tetrahedrons have some little cubes in common.

It follows one example of a cube divided in tetrahedrons without overlaps.

**Example 4'.**

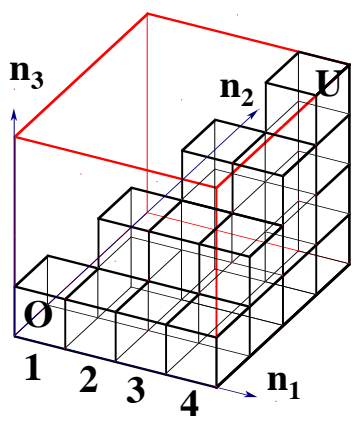

$C^3_4 = \{(X_1,X_2,X_3)$: $X_h = 1,2,3,4$, $h=1,2,3\} = 4^3$, contains the tetrahedrons characterized from the 3! permutations, $X_{h_1}\, X_{h_2}\, X_{h_3}$, of $X_1\, X_2\, X_3$, to which the *rule of the climbs* is applied.

We call *factorial tetrahedrons* the tetrahedrons[7] previously characterized. They are the following (where, for brevity, as an example, the little cube (1,2,1) will be indicated with 121, and the ends of the fishbone are omitted)

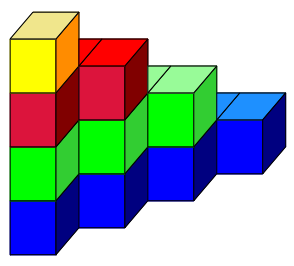

- 1) one tetrahedron $T^3_4 = \{(X_1,X_2,X_3)$: $X_1 \geq X_2 \geq X_3\} = \frac{4^{3\geq}}{3!} = \frac{4 \cdot 5 \cdot 6}{3!} =$

$= \{(111\ 211\ 311\ 411\ 221\ 321\ 421\ 331\ 431\ 441)\ (222\ 322\ 422\ 332\ 432\ 442)\ (333\ 433\ 443)\ (444)\} =$

$= T^2_4 \cup T^2_3 \cup T^2_2 \cup T^2_1$, where $T^2_h$ is on the plane: $X_3 = 5-h$, h=1,2,3,4.

- 4 tetrahedrons $T^3_3 = T^2_3 \cup T^2_2 \cup T^2_1 = \frac{3^{3\geq}}{3!} = \frac{3 \cdot 4 \cdot 5}{3!}$

2) $X_1 \geq X_3 > X_2$: {(212 213 313 214 314 414) (323 324 424) (434)}, (on the planes $X_2$=1, 2, 3)

3) $X_2 > X_1 \geq X_3$: {(121 131 141 231 241 341) (242 342 442) (443)}, (on the planes $X_3$=1, 2, 3)

4) $X_2 \geq X_3 > X_1$: {(122 123 124 133 134 144) (233 234 244) (344)}, (on the planes $X_1$=1, 2, 3)

5) $X_3 > X_1 \geq X_2$: {(112 113 213 114 214 314) (223 224 324) (334)}, (on the planes $X_2$=1, 2, 3)

---

[6] What it has been asserted can be generalized to any tetrahedron transforming by an affinity a cube and the contained tetrahedrons in a parallelepiped, knowing its measure.

[7] If it is important to distinguish the factorial tetrahedrons from other types of tetrahedrons, $T^d_n$, can be used the symbol $T!^d_n$.

- and one tetrahedron $T^3{}_2 = \frac{2^{3\geq}}{3!} = \frac{2\cdot 3\cdot 4}{3!}$

6) $X_3>X_2>X_1$: {(123 124 134) (234)}, (on the planes $X_1$=1, 2).

In this way we can verify that the previous fishbones represent an Euler's Escher that operates a partition of $C^3{}_4$.

## 5. OTHER TYPES OF TETRAHEDRONS IN THE SPACE $\underline{N}^d$

**Lemma 4**. *Each little cube of the fishbone*

$\mathbf{I}=\{(X_1,X_2,...,X_d)$: $\mathbf{b}\geq X_1\geq...\geq X_d\geq\mathbf{a+1}$, $X_h$=a+1,a+2,..., h=1,2,...,d}, a,b$\in\mathbf{N}$

*correspond to a little cubes of the fishbone.* **I** a **I'** have the same number of solutions.

$\mathbf{I'}=\{(X_1,X_2,...,X_d)$: $\mathbf{b-a}\geq X_1\geq...\geq X_d\geq\mathbf{1}$, $X_h$=**1**,2,..., h=1,2,...,d}.

*Proof.* Each little cube in **I'** and in **I** is obtained by the translation respectively defined subtracting or adding the number **a** from the coordinates of each little cube of **I** and **I'**.

**Lemma 5**. *The set of little cubes of the fishbone*

$\mathbf{J}=\{(X_1,X_2,...,X_d)$: $\mathbf{n}\geq X_1\geq...\geq X_d\geq\mathbf{1}$, $X_h$=1,2,..., h=1,2,...,d}

*represents a tetrahedron in d dimensions with edges of length n,* $\mathbf{T^d{}_n}$.

*Proof.* $\mathbf{n}\geq X_1\geq X_2\geq\mathbf{1}$ is a tetrahedron in 2 dimensions, $T^2{}_n$, $\forall$n.

Hip.: supposed that $\mathbf{h}\geq X_1\geq...\geq X_{d-1}\geq\mathbf{1}$ is a tetrahedron in d-1 dimensions, $\mathbf{T^{d-1}{}_h}$, $\forall$h, applying lemma 4, it turns out

$$(\mathbf{n}\geq X_1\geq...\geq X_d\geq\mathbf{1}) = \bigcup_{h=1}^{n} (\mathbf{n}\geq X_1\geq...\geq X_d=\mathbf{h}) = \bigcup_{h=1}^{n} (\mathbf{n}\geq X_1\geq...\geq X_{d-1}\geq\mathbf{h}) = \bigcup_{1}^{n}\mathbf{T^{d-1}{}_k}, \; k=n-h.$$

The set **J** represents a tetrahedron in d dimensions because it is valid the *rule of the sections*, and since $\mathbf{n}\geq X_h\geq\mathbf{1}$, $\forall$h, than the length of its edges is n. Than $\mathbf{J} = \mathbf{T^d{}_n} = \frac{\mathbf{n^{d\geq}}}{\mathbf{d!}}$,[8]

Example $5\geq X_1\geq X_2\geq X_3\geq 1$ 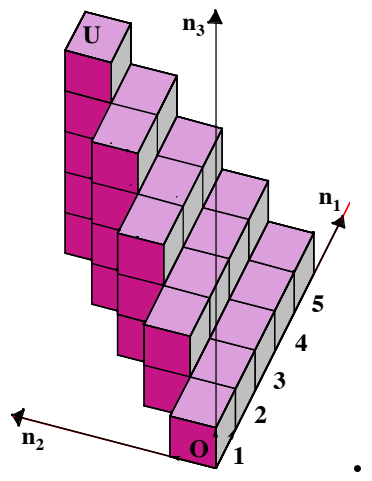.

**Lemma 5a.** *The tetrahedron* $\{(X_1,X_2,...,X_d)$: $\mathbf{n}\geq X_1\geq...\geq X_d\geq\mathbf{1}$, $X_h$=**1,2,**..., h=1,2,...,d} *has the same sections of the rectangular tetrahedron* $\{(X_1,X_2,...,X_d)$: $\sum_1^d \mathbf{X_j}\leq\mathbf{n}$, $X_j$=**1,2,**... j=1,2,...,d}.

Example $\sum_1^3 \mathbf{x_j}\leq\mathbf{4}$ 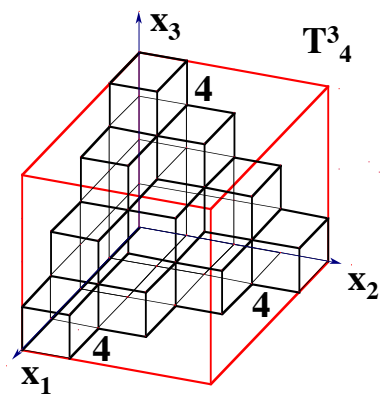

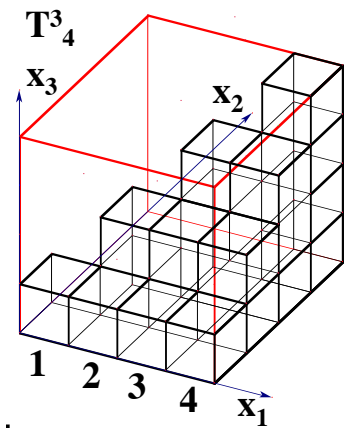

$4\geq X_1\geq X_2\geq X_3\geq 1$

*Proof.* The lemma is true considering a single coordinate $\{(X_1)$: $\mathbf{n}\geq X_1\geq\mathbf{1}$, $X_1$=**1,2,**...,}, $\forall$n.

Hip.: supposed that lemma 5a is true in d-1 dimensions, $\forall$n, this lemma is also valid in d dimensions, applying lemma 5, $(\mathbf{n}\geq X_1\geq...\geq X_d\geq\mathbf{1}) = \bigcup_1^n \mathbf{T^{d-1}{}_k}$.

We introduce what follows with two examples.

---

[8] It follow that, dividing a cube $C^d{}_e$, $\mathbf{e}\in$R, in $n^d$ little cubes, $\mathbf{T^d{}_n} = \frac{\mathbf{n^{d\geq}}}{\mathbf{d!}}\frac{\mathbf{e^d}}{\mathbf{n^d}}$ tends to $\frac{\mathbf{e^d}}{\mathbf{d!}}$, and this method is generalized, with an affinity, to a parallelepiped, knowing his volume.

**Example 5**.
The little cubes of the fishbone $\mathbf{6}\ge X_2 > X_1 \ge X_3 \ge \mathbf{1}$ indicate a tetrahedron that is a translation of the tetrahedron $T^3{}_5 = \mathbf{5}\ge X_2 \ge X_1 \ge X_3 \ge \mathbf{1}$.
The solutions of the fishbone ($= \bigcup_{1}^{5}(\mathbf{6 \ge X_2 > X_1 \ge X_3 = h})$) are

121 131 141 151 161
231 241 251 261    232 242 252 262
341 351 361    342 252 362    343 353 363
451 461    452 462    453 463    454 464
561    562    563    564    565

the solutions ($= \bigcup_{1}^{5}(\mathbf{5 \ge X_2 \ge X_1 \ge X_3 = h})$) of the tetrahedron $T^3{}_5$ ($= \bigcup_{1}^{5}(\mathbf{5 \ge X_2 \ge X_1 \ge X_3 = h})$) are

111 121 131 141 151
221 231 241 251    222 232 242 252
331 341 351    332 342 352    333 342 353
441 451    442 452    443 453    444 454
551    552    553    554    555

There is a translation because the solutions of $\mathbf{6}\ge X_2 > X_1 \ge X_3 \ge \mathbf{1}$ can be obtained from the solutions of $T^3{}_5$, transforming $X_2$ in $X_2+1$ and leaving the other coordinates unchanged.

**Example 6**.
The fishbone **H**: $\mathbf{4}\ge X_1\ge X_3>X_2\ge X_4\ge X_8>X_5\ge X_7>X_6\ge\mathbf{1}$, has only the solution **43432123**, that correspond to the solution of the fishbone obtained subtracting 3 from his extreme 4, **K**: $\mathbf{1}\ge X_1\ge X_3\ge X_2\ge X_4\ge X_8\ge X_5\ge X_7\ge X_6\ge\mathbf{1}$, where 3 is the number of the signs " > " between the variables in H, replaced with the signs " ≥" in K.

The only one solution of K, 1111 1111, that indicates a tetrahedron $T^8{}_1$, is obtained transforming: $X_6$, $X_7$, $X_5$, $X_8$, $X_4$, $X_2$, $X_3$, $X_1$, in H, respectively in: $X_6$, $X_7$-1, $X_5$-1, $X_8$-2, $X_4$-2, $X_2$-2, $X_3$-3, $X_1$-3, in K.

In other words, we can say that the solutions of the previous fishbone

**H**: $\mathbf{4}\ge X_1\ge X_3 > X_2\ge X_4\ge X_8 > X_5\ge X_7 > X_6\ge\mathbf{1}$ correspond to the solutions of
$\mathbf{4}\ge X_1\ge X_3 > X_2\ge X_4\ge X_8 > X_5\ge X_7 \ge X_6\ge\mathbf{2}$, that correspond to the solutions of
$\mathbf{4}\ge X_1\ge X_3 > X_2\ge X_4\ge X_8 \ge X_5\ge X_7 \ge X_6\ge\mathbf{3}$, that correspond to the solutions of
$\mathbf{4}\ge X_1\ge X_3 \ge X_2\ge X_4\ge X_8 \ge X_5\ge X_7 \ge X_6\ge\mathbf{4}$, that, applying lemma 4, correspond to the solutions of
**K**: $\mathbf{1}\ge X_1\ge X_3 \ge X_2\ge X_4\ge X_8 \ge X_5 \ge X_7 \ge X_6\ge\mathbf{1}$.

Repeating in general terms what it is now indicated, we have the following lemma.

**Lemma 6.**
The solutions of the fishbone **V**: $\mathbf{n}\ge X_{h_1}\ge X_{h_2}> X_{h_3}> X_{h_4}\ge ...> X_{h_{d-1}}> X_{h_d}\ge\mathbf{1}$,
*in which there are* ***v*** *signs* ">", *is one-to-one correspondent to the solutions of the tetrahedron* **Z**: $\mathbf{n\text{-}v} \ge X_{h_1}\ge X_{h_2}\ge X_{h_3}\ge X_{h_4}\ge ...\ge X_{h_{k-1}}\ge X_{h_d}\ge\mathbf{1}$, *where each of the* ***v*** *signs* ">" *of the fishbone* **V** *have been replaced in* **Z** *with a sign* "≥".

*Proof.* Since for $X_h$, and $\mathbf{c}\in N$, $X_h>c$ is equivalent to $X_h\ge c+1$, $\forall \mathbf{h}$, proceeding from right to the left in the fishbone **V**, and removing an unit from the right extreme for each sign " > " that it is met, the following fishbone one-to-one correspondent to the fishbone **V**, are respectively determined $\mathbf{n}\ge X_{h_1}\ge X_{h_2}> X_{h_3}> X_{h_4}\ge ...> X_{h_{k-1}}\ge X_{h_d}\ge\mathbf{2}$,
$\mathbf{n}\ge X_{h_1}\ge X_{h_2}> X_{h_3}> X_{h_4}\ge ...\ge X_{h_{k-1}}\ge X_{h_d}\ge\mathbf{3}$, ... $\mathbf{n}\ge X_{h_1}\ge X_{h_2}\ge X_{h_3}\ge X_{h_4}\ge ...\ge X_{h_{k-1}}\ge X_{h_d}\ge\mathbf{1+v}$
that is, applying lemma 4, one-to-one correspondent to the tetrahedron **Z**.

**Definition**. We call **factorial tetrahedrons**, the set of little cubes of this kind of fishbone $\mathbf{n} \geq x_{h_1} \geq x_{h_2} > x_{h_3} > x_{h_4} \geq ... \geq x_{h_d} \geq \mathbf{1}$, where $x_{h_1}, x_{h_2}, x_{h_3}, ..., x_{h_d}$ is a permutation of variables and the *rule of the climbs* is valid.

**Theorem.** *We have* $C^d{}_n = n^d = \sum_{s=0}^{d-1} \left\langle {d \atop s} \right\rangle T^d_{n-s}$ $\forall n; d \in \mathbf{N}$, and, to add powers of integers, 1, 2 or t times, we have $\sum_{j=1}^{n} j^d = \sum_{s=0}^{d-1} \left\langle {d \atop s} \right\rangle T^{d+1}_{n-s}$, $\sum_{m=1}^{n} \sum_{j=1}^{m} j^d = {}^2\sum_{j=1}^{n} j^d = \sum_{s=0}^{d-1} \left\langle {d \atop s} \right\rangle T^{d+2}_{n-s}$, ${}^t\sum_{j=1}^{n} j^d = \sum_{s=0}^{d-1} \left\langle {d \atop s} \right\rangle T^{d+t}_{n-s}$

*where Eulerian number*[9] $\left\langle {d \atop s} \right\rangle$ *is the number of the factorial d-dimensional tetrahedrons with edges of length* n-s, $\mathbf{T}^d{}_{\mathbf{n-s}}$, *contained in* $C^d{}_n$.

*Proof.* All the d! *factorial tetrahedrons* $\mathbf{n} \geq x_{h_1} \geq x_{h_2} > x_{h_3} > x_{h_4} \geq ... \geq x_{h_d} \geq \mathbf{1}$, represent an Escher and therefore they decompose the cube $C^d{}_n$. Considering that in the permutations of 1,2,3,…,d, the number $\left\langle {d \atop s} \right\rangle$ of those with **s** signs "≥" is equal to the number of the permutations with **s** signs ">", it turns out $C^d{}_n = \mathbf{n}^\mathbf{d} = \sum_{s=0}^{d-1} \left\langle {d \atop s} \right\rangle T^d_{n-s}$, where $\left\langle {d \atop s} \right\rangle$, considering the previous lemmas, becomes the number of the factorial d-dimensional tetrahedrons with edges of length n-s, $\mathbf{T}^d{}_{\mathbf{n-s}}$, contained in $C^d{}_n$, expressed by the fishbones with **s** signs ">".

Applying the property $T^d{}_n = \sum_{1}^{n} T^{d-1}{}_h$, 1, 2 or t times, to $\mathbf{n}^\mathbf{d} = \sum_{s=0}^{d-1} \left\langle {d \atop s} \right\rangle T^d_{n-s}$, we have

$$\sum_{j=1}^{n} j^d = \sum_{s=0}^{d-1} \left\langle {d \atop s} \right\rangle T^{d+1}_{n-s}, \quad \sum_{m=1}^{n} \sum_{j=1}^{m} j^d = \sum_{s=0}^{d-1} \left\langle {d \atop s} \right\rangle T^{d+2}_{n-s} = {}^2\sum_{j=1}^{n} j^d, \quad {}^t\sum_{j=1}^{n} j^d = \sum_{s=0}^{d-1} \left\langle {d \atop s} \right\rangle T^{d+t}_{n-s}$$

- *A not geometric synthetic proof of* $\mathbf{n}^\mathbf{d} = \sum_{s=0}^{d-1} \left\langle {d \atop s} \right\rangle T^d_{n-s}$.

The Worpitzky's identity (1883) $\mathbf{x}^\mathbf{d} = \sum_{k=0}^{d-1} \left\langle {d \atop k} \right\rangle \binom{x+k}{d}$, $d \in \mathbf{N}$, since $\left\langle {d \atop k} \right\rangle = \left\langle {d \atop d-k-1} \right\rangle$, integer k>0, and putting s=d-k-1 and x= $n \in \mathbf{N}$, becomes $\mathbf{n}^\mathbf{d} = \sum_{s=0}^{d-1} \left\langle {d \atop s} \right\rangle \binom{n-s+d-1}{d} = \sum_{s=0}^{d-1} \left\langle {d \atop s} \right\rangle T^d_{n-s}$.

## 7. EULER TRIANGLE AND EXAMPLES

| d | $\left\langle {d \atop 0} \right\rangle$ | $\left\langle {d \atop 1} \right\rangle$ | $\left\langle {d \atop 2} \right\rangle$ | $\left\langle {d \atop 3} \right\rangle$ | $\left\langle {d \atop 4} \right\rangle$ | $\left\langle {d \atop 5} \right\rangle$ | $\left\langle {d \atop 6} \right\rangle$ | $\left\langle {d \atop 7} \right\rangle$ | $\left\langle {d \atop 8} \right\rangle$ … |
|---|---|---|---|---|---|---|---|---|---|
| 0 | 1 | 0 | ... | | | | | | |
| 1 | 1 | 0 | ... | | | | | | |
| 2 | 1 | 1 | 0 | ... | | | | | |
| 3 | 1 | 4 | 1 | 0 | ... | | | | |
| 4 | 1 | 11 | 11 | 1 | 0 | ... | | | |
| 5 | 1 | 26 | 66 | 26 | 1 | 0 | ... | | |
| 6 | 1 | 57 | 302 | 302 | 57 | 1 | 0 | ... | |
| 7 | 1 | 120 | 1191 | 2416 | 1191 | 120 | 1 | 0 | |
| 8 | 1 | 247 | 4293 | 15619 | 15619 | 4293 | 247 | 1 | 0 |
| 9 | 1 | 502 | 14608 | 88234 | 156190 | 88234 | 14608 | 502 | 1 |
| 10 | 1 | 1013 | 47840 | 455192 | 1310354 | 1310354 | 455192 | 47840 | 1013... |
| 11 | 1 | 2036 | 152637 | 2203488 | 9738114 | 15724248 | 9738114 | 2203488 | 152637... |
| 12 | 1 | 4083 | 478271 | 10187685 | 66318474 | 162512286 | 162512286 | 66318474 | 10187685... |
| 13 | 1 | 8178 | 1479726 | 45533450 | 423281535 | 1505621508 | 2275172004 | 1505621508 | 423281535... |
| 14 | 1 | 16369 | 4537314 | 198410786 | 2571742175 | 12843262863 | 27971176092 | 27971176092 | 12843262863... |
| 15 | 1 | 32752 | 13824739 | 848090912 | 15041229521 | 102776998928 | 311387598411 | 447538817472 | 311387598411... |
| 16 | 1 | 65519 | 41932745 | 3572085255 | 85383238549 | 782115518299 | 3207483178157 | 6382798925475 | 6382798925475... |
| 17 | 1 | 131054 | 126781020 | 14875399450 | 473353301060 | 5717291972382 | 31055652948388 | 83137223185370 | 114890380658550... |

[9] The Eulerian numbers: $\left\langle {d \atop s} \right\rangle$, d≥0, are defined by: $\left\langle {d \atop s} \right\rangle = (1+s) \left\langle {d-1 \atop s} \right\rangle + (d-s) \left\langle {d-1 \atop s-1} \right\rangle$, $\left\langle {1 \atop 0} \right\rangle = 1$, $\left\langle {d \atop k} \right\rangle = 0$, k≥d.

**Examples 7.**

$4^3 = T^3{}_4 + 4\,T^3{}_3 + T^3{}_2 = 20 + 40 + 4 = 64$

$1^4 = T^4{}_1 = 1,\ 2^4 = T^4{}_2 + 11T^4{}_1 = 5+11=16$
$3^4 = T^4{}_3 + 11T^4{}_2 + 11T^4{}_1 = 15+55+11=81$
$4^4 = T^4{}_4 + 11T^4{}_3 + 11T^4{}_2 + T^4{}_1 = 35 + 165 + 55+1 = 256$
$5^4 = T^4{}_5 + 11T^4{}_4 + 11T^4{}_3 + T^4{}_2 = 70+385+165+5 = 625$

- $\sum_{j=1}^{5} j^4 = 1 + 2^4 + 3^4 + 4^4 + \mathbf{5}^4 = 1 + 16 + 81 + 256 + 625 = 979 =$
  $= T^5{}_5 + 11T^5{}_4 + 11T^5{}_3 + T^5{}_2 = 126+616+231+6=979$
- $\sum_{m=1}^{5}\sum_{j=1}^{m} j^4 = 1+(1+16)+(1+16+81)+(1+16+81+256)+(1+16+81+256+625)=1449=$
  $= T^6{}_5 + 11T^6{}_4 + 11T^6{}_3 + T^6{}_2 = 210+924+308+7=1449$
- $\sum_{j=1}^{100} j^4 = T^5{}_{100} + 11T^5{}_{99} + 11T^5{}_{98} + T^5{}_{97} =$

$= \frac{1}{5!}(\mathbf{100}\cdot\mathbf{101}\cdot 102\cdot 103\cdot 104 + 11\cdot 99\cdot\mathbf{100}\cdot\mathbf{101}\cdot 102\cdot 103 + 11\cdot 98\cdot 99\cdot\mathbf{100}\cdot\mathbf{101}\cdot 102 + 97\cdot 98\cdot 99\cdot\mathbf{100}\cdot\mathbf{101})=$

$= \frac{100\cdot 101}{5!}(102\cdot 103\cdot 104+11\cdot 99\cdot 102\cdot 103+11\cdot 98\cdot 99\cdot 102+97\cdot 98\cdot 99) =$

$= \frac{100\cdot 101}{5!}\ 24\ 360\ 396 = 2\ 050\ 333\ 330 = 2.5.41.67.101.739.$

**Theorem 2.** *It turns out* $\sum_{j=1}^{n} j^d = \frac{n(n+1)}{(d+1)!}k,\ k\in\mathbf{N},$

*Proof.* It is sufficient to generalize the previous example.

## 8. OTHER PRESENCES OF THE EULERIAN NUMBERS IN THE CUBES

**Theorem 3.** *We have* $\left\langle {d \atop s} \right\rangle = \nabla^{d+1}(s+1)^d = (s+1)^d - \binom{d+1}{1}s^d + \binom{d+1}{2}(s-1)^d \ldots$ [10]

*Introduction.* We indicate with $R^d(0)$ the part of the space $\mathbf{R}^d$, with the points with not negative coordinates. The vertex of $R^d(0)$ is the origin, **O**.

We consider the $H=\binom{d}{h}$ vertices of $C^d(e)$, $V_h$, with h coordinates equal to e and the others null, with an order $V_{h_1}, V_{h_2}, \ldots, V_{h_j}, \ldots, V_{h_H}$, $\forall h=1,2,\ldots,d$. At the end we place $R^d(h_j)$ the part of the space obtained translating $R^d(0)$ in $V_{h_j}$ with the vector $OV_{h_j}$.

**Lemma 7.** *The measure* $T^d(\le s) \cap R^d(h_j)$ *is equal to that of* $T^d(\le s\text{-}he)$, $\forall s\in R$, $\forall h, j\in N$, where $T^d(\le s)=\{(x_1,x_2,\ldots,x_d): \sum_1^d x_j \le s,\ 0\le x_h,\ h=1,2,\ldots,d\}=\frac{s^d}{d!}$.[11]

*Proof.* We have to consider in $T^d(\le s) \cap R^d(h_j)$ the h coordinates greater or equal to e of $V_{h_j}$. Taken one of these, as an example $x_i$, we can replace $x_i$+e to $x_i$ in the definition of $T^d(\le s)$, so $\sum_1^d x_j \le s$ becomes $\sum_1^d x_j \le s\text{-}e$, and $\sum_1^d x_j \le s\text{-}he$, considering all the h coordinates greater or equal to **e** of $V_{h_j}$.

---

[10] Where $\nabla f(x)= f(x)-f(x-1)$, and $\nabla^h f(x)= f(x)-\binom{h}{1} f(x-1)+\binom{h}{2} f(x-2) - \ldots$

[11] It's very easy to prove that $T!^d(s)$ and $T^d(\le s)$ have the same measures.

**Lemma 8**. $C^d(e;\le s) = C^d(e) \cap T^d(\le s) =$

$$= T^d(\le s) - \binom{d}{1} T^d(\le s\text{-}e) + \binom{d}{2} T^d(\le s\text{-}2e) - \binom{d}{3} T^d(\le s\text{-}3e) + \ldots = \nabla_e^{\,d} T^d(\le s) = T^d(\le s|e),$$

*where* $T^d(\le r) = 0$, *when* $r \le 0$, $C^d(e)=\{(x_1,x_2,\ldots,x_d): 0\le x_h \le e,\ h=1,2,\ldots,d\}$, $e\in \mathbf{R}$, $s\in \mathbf{R}$.

**Example 8**. $C^2(e;\le s) = T^2(\le s) - 2T^2(\le s\text{-}e) + T^2(\le s\text{-}2e) = \nabla_e^{\,2} T^2(\le s) = T^2(\le s|e)$.

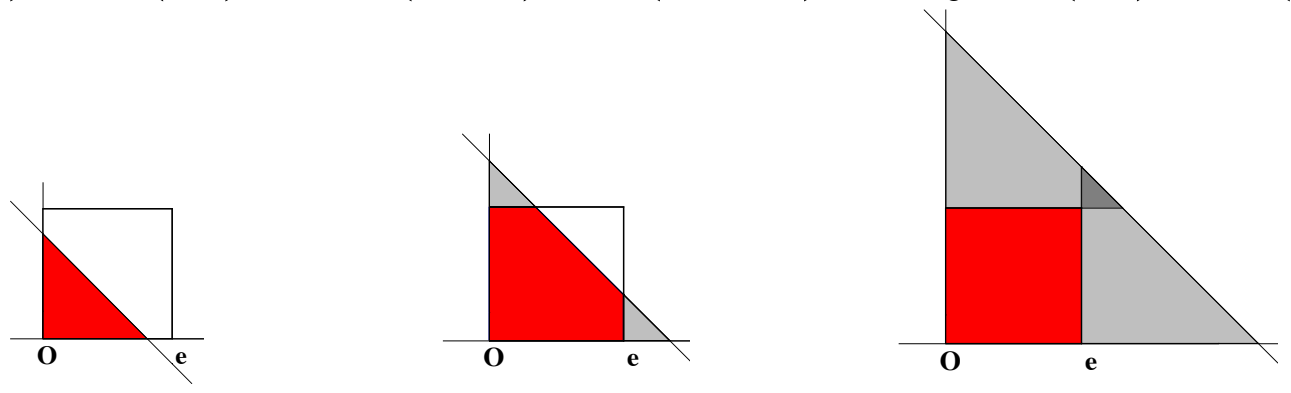

**Example 9**. $C^3(e;\le s) = T^3(\le s|e) = T^3(\le s) - 3T^3(\le s\text{-}e) + 3T^3(\le s\text{-}2e) - T^3(\le s\text{-}3e)$, $\forall s \in R$,

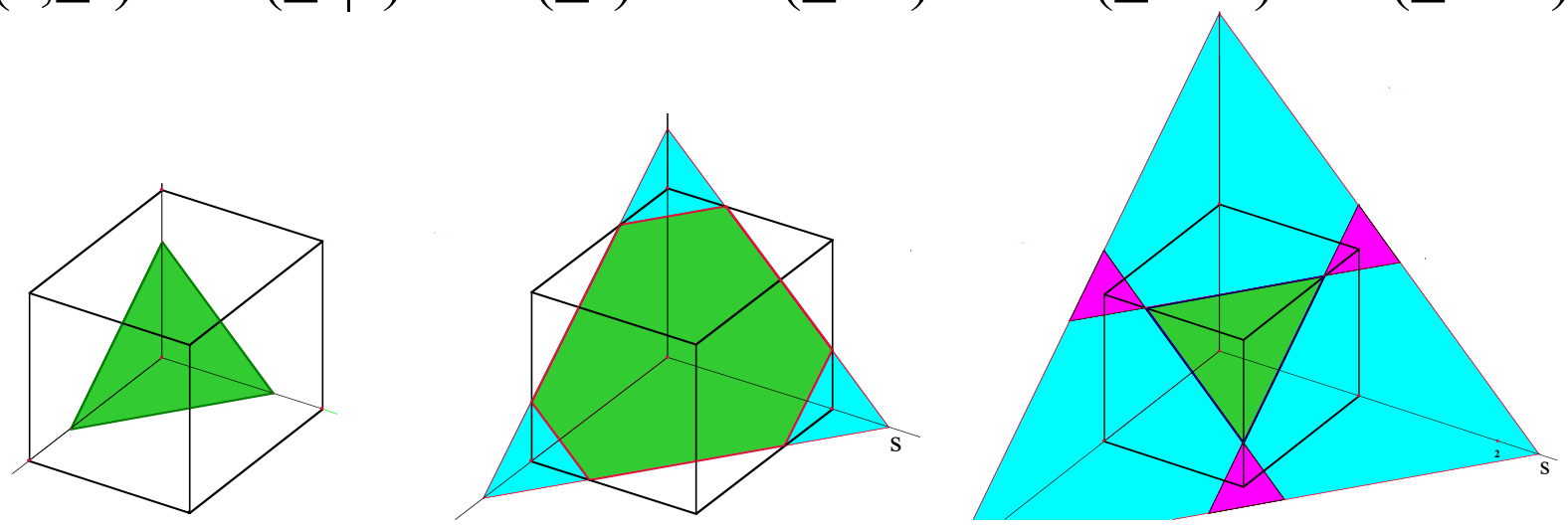

*Proof*. We consider the events[12] $E_h=(0\le x_h\le e) = (0\le x_h)-(x_h>e)$, $h=1,2,\ldots,d$, and synthetically indicate with $C^d(e)$ the event "the point $P\in C^d(e)$". It turns out that

$C^d(e) = (0\le x_1\le e)\cap \ldots \cap (0\le x_d\le e)= E_1E_2\ldots E_d =$

$$= ((0\le x_1)-(x_1>e))\cdot((0\le x_2)-(x_2>e))\cdot\ldots\cdot((0\le x_d)-(x_d>e)) =$$

$$= R^d(0)-[\,R^d(1_1)+R^d(1_2)+\ldots]+[\,R^d(2_1)+R^d(2_2)+\ldots]+\ldots$$

$$+(-1)^h[\,R^d(h_1)+R^d(h_2)+\ldots] + \ldots$$

Considering $C^d(e;\le s) = C^d(e)\cap T^d(\le s)$ and replacing $T^d(\le s\text{-}he)$ to $T^d(\le s)\cap R^d(h_j)$ (lemma 7), remembering that the number of $R^d(h_j)$ is $\binom{d}{h}$, and $T^d(\le s)=0$, $s\le 0$, we have

$$C^d(e;\le s) = C^d(e)\cap T^d(\le s) = T^d(\le s) - \binom{d}{1}T^d(\le s\text{-}e) + \binom{d}{2} T^d(\le s\text{-}2e) + \ldots(-1)^d\binom{d}{d} T^d(\le s\text{-}de) =$$

$$= \nabla_e^{\,d}\, T^d(\le s) = T^d(\le s|e).$$ [13]

**Corollary**. *It turns out* $C^d(e;s\text{-}e\le x\le s)= C^d(e;\le s) - C^d(e;\le s\text{-}e)= \nabla_e^{\,d+1}\dfrac{s^d}{d!}$.

*Proof.* $C^d(e;s\text{-}e\le x\le s) = \nabla_e C^d(e;\le s) = \nabla_e^{\,d} T^d(\le s) - \nabla_e^{\,d} T^d(\le s\text{-}e) = \nabla_e^{\,d+1} T^d(\le s) = \nabla_e^{\,d+1}\dfrac{s^d}{d!}$.

In particular, $\forall s\in \mathbf{N}$, we obtain the volumes of the portions of $\mathbf{C^d(e)}$ between the sections on its vertices $\sum_1^d x_j = (s-1)e$, $\sum_1^d x_j = se$.

We consider to simplify $e = 1$ and indicate $C^d(1;s\text{-}1\le x\le s)$ with $C^d(s\text{-}1\le x\le s)$ and $\nabla_1$ with $\nabla$.

*Proof of the Theorem 3.* We have from the previous corollary and for $s\in \mathbf{N}$

$$d!C^d(s\text{-}1\le x\le s) = \nabla^{d+1} s^d = \sum_0^{d+1}(-1)^h\binom{d+1}{h}(s-h)^d = \sum_0^{d+1}(-1)^h\frac{d!(s-h)^{d-1}}{h!(d+1-h)!}(d+1)(s-h) = \gamma)$$

---

[12] If A and B are events, B-A is an event if and only if $A\subset B$.

[13] The proof with $s\in R$, is also a way to calculate the density and the probability distribution of the sum of independent random numbers, identically distributed, with uniform distribution (Barra, 2008).

Placing $\underline{(d+1)(s-h)}$=ds-dh+s-h-sh+sh=$\underline{s(d+1-h)-h(d+1-s)}$, and considering that d+1-h = 0 when h = d+1, and that h(d+1-s) = 0 for h = 0, we obtain

$$\gamma) = s\sum_{h=0}^{d}(-1)^h\binom{d}{h}(s-h)^{d-1} - (d+1-s)\sum_{h=1}^{d+1}(-1)^h\frac{d!(s-h)^{d-1}}{(h-1)!(d+1-h)!},$$

from which, placing k+1 = h, it is obtained

$$\gamma) = \nabla^{d+1}s^d = s\sum_{h=0}^{d}(-1)^h\binom{d}{h}(s-h)^{d-1} + (d+1-s)\sum_{k=0}^{d}(-1)^h\frac{d!(s-k-1)^{d-1}}{(k)!(d-k)!} =$$

$= s\,\nabla^d s^{d-1} + (d+1-s)\ \nabla^d (s-1)^{d-1}$. Therefore $\nabla^{d+1}s^d$ verifies the recurrence of $\left\langle {d \atop s-1} \right\rangle$ and their initial conditions. So we have $\left\langle {d \atop s-1} \right\rangle = \nabla^{d+1}s^d$, or $\left\langle {d \atop s} \right\rangle = \nabla^{d+1}(s+1)^d = \delta)$.

**Example 10**. $\left\langle {d \atop 0} \right\rangle$=**1**,

$\left\langle {2 \atop 0} \right\rangle$=**1**, $\left\langle {2 \atop 1} \right\rangle=2^2-3\cdot1=$**1**, $\left\langle {2 \atop 2} \right\rangle=3^2-3\cdot2^2+3\cdot1=$ **0**

$\left\langle {3 \atop 0} \right\rangle$=**1**, $\left\langle {3 \atop 1} \right\rangle=2^3-4\cdot1=$**4**, $\left\langle {3 \atop 2} \right\rangle=3^3-4\cdot2^3+6\cdot1=$**1** $\left\langle {3 \atop 3} \right\rangle=4^3-4\cdot3^3+6\cdot2^3-4\cdot1=$ **0**

$\left\langle {4 \atop 0} \right\rangle$=**1**, $\left\langle {4 \atop 1} \right\rangle=2^4-5\cdot1=$**11**, $\left\langle {4 \atop 2} \right\rangle=3^4-5\cdot2^4+10\cdot1=$**11** $\left\langle {4 \atop 3} \right\rangle=4^4-4\cdot3^4+6\cdot2^4-4\cdot1=$ **0**

We have theorem 3 considering various cubes with increasing vertices. Considering only one cube, we have the following corollary.

**Corollary 2**. *We find Eulerian numbers in the ratio between the volumes of the portions of* $\mathbf{C^d(e)}$ *placed between the sections on its vertices* $\sum_1^d x_j = (s-1)e$, $\sum_1^d x_j = se$, s=1,2,…,d, *or, evidently, in the analogous portions of* $\mathbf{C^d(\sqrt[d]{d!}) = d!}$, *or in the ratio between the measures of the previous sections on the vertices of* $\mathbf{C^d(e)}$.

*Proof.* From **δ**), dividing and multiplying by d!, we have $\left\langle {d \atop s} \right\rangle$=d! $\nabla^{d+1}T^d(\le s+1)$, and since $T^d(\le s+1) = \frac{T^d(s+1)}{\sqrt{d+1}}$,[14] we obtain $\left\langle {d \atop s} \right\rangle$ = d! $\nabla^{d+1}\frac{T^d(s+1)}{\sqrt{d+1}}$, and dividing by $T^d(1)=\frac{\sqrt{d+1}}{d!}$, we also find Eulerian numbers in the ratio between the measures of sections on vertexes of $C^d_e$.

**Example 11**. *Measures of portions*

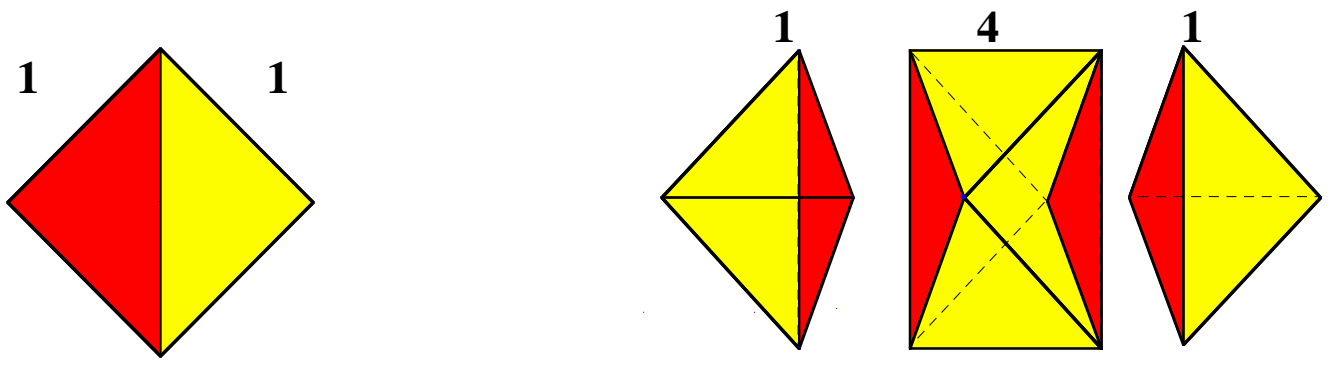

**Example 12**. *Measures of sections*

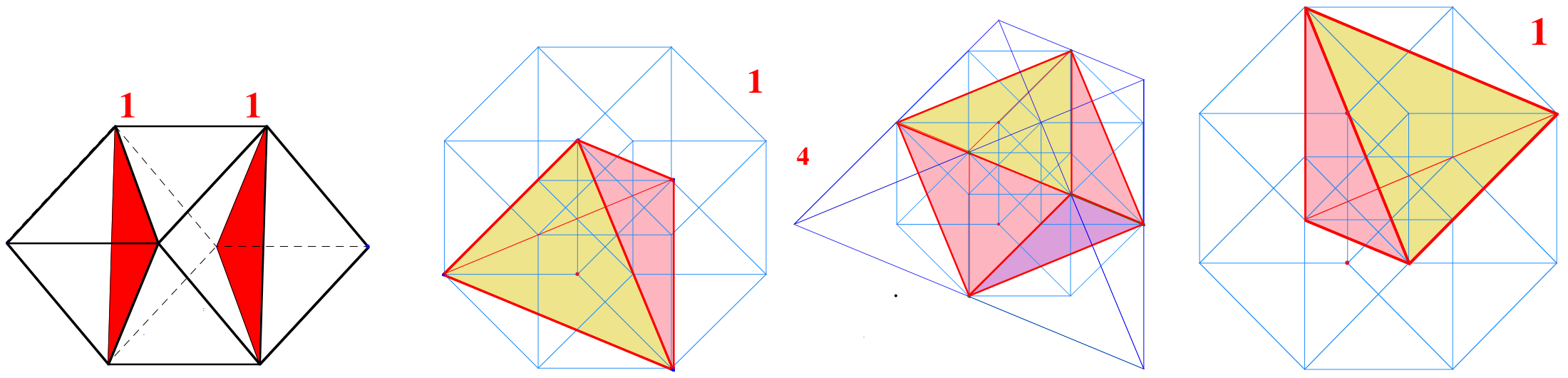

[14] It's very easy to prove that $T^d(s)=\sqrt{d+1}\ T^d(\le s)$.

## 9. SHORT COMPARISON WITH OTHER METHODS

The most famous numbers used to sum powers are the Bernoulli (1654-1705) numbers $B_0, B_1, B_2, B_3, B_4, \ldots$, where $B_n = 0$ for all odd $\boldsymbol{n}$ other than 1. $B_1 = ½$ or $-1/2$ depending on the convention adopted. The values of the first few nonzero Bernoulli numbers are

| n | 0 | 1 | 2 | 4 | 6 | 8 | 10 | 12 |
|---|---|---|---|---|---|---|---|---|
| $B_n$ | 1 | ±1/2 | 1/6 | −1/30 | 1/42 | −1/30 | 5/66 | −691/2730 |

The Bernoulli numbers appear in the Taylor series expansions of the tangent and hyperbolic tangent functions, in the Euler–Maclaurin formula, and in expressions for certain values of the Riemann zeta function ... , but to have a new $B_k$ is necessary to make reference to all the $B_h$, with h<k. Thus a "triangle" with Bernoulli numbers is not easy to build.

We think that the sum of powers obtained through the Euler's triangle in a geometric way is simple and original.

Moreover the numbers of Euler are integers and have many properties, easy to prove.

## BIBLIOGRAPHY

Barra M., 2001, Ipersolidi e altri argomenti, *Progetto Alice,* N. 5, Vol. 2, Ed. Pagine, 191-246.

Barra M., 2008, Matematica e software di geometria dinamica seguendo le indicazioni scientifiche e didattiche di Bruno de Finetti, *Progetto Alice,* N. 26, Vol. 9, Ed. Pagine, 191-230.

Barra M., 2009, Numeri Euleriani, Box-numbers, Numeri di Stirling di seconda specie e Triangolo Aritmetico di Pascal per ottenere le somme delle potenze degli interi. Approfondimenti e collegamenti, *Progetto Alice,* N. 28, Vol. 10, Ed. Pagine, 97-115.

Comtet, L., 1974, Permutations by Number of Rises; Eulerian Numbers. §6.5 in *Advanced Combinatorics: The Art of Finite and Infinite Expansions, rev. enl. ed.* Dordrecht, Netherlands: Reidel, 240-246.

Coxeter H.S.M., 1963 (1948), *Regular polytopes*, The Macmillan Company, New York.

De Fermat P., 1959 (1670), *Osservazioni su Diofanto*, Boringhieri.

Edwards A.W.F., 2002 (1987), *Pascal's Arithmetical Triangle*, Charles Griffin & Company Limited, London.

Edwards A.W.F, 1982, Sums of powers of integers-a little of the history, *The Mathematical Gazette*, 66, 22-28.

Edwards A.W.F, 1986, A quick route to sums of powers, *American Mathematical Monthly*, n. 93, 451-455.

Graham R.L., Knuth D.E., Patashnik O., 1989, *Concrete Mathematics.* Addison-Wesley, Reading, MA.

Worpitzky J., 1883, Studien über die Bernoullischen und Eulerschen Zahlen, *Journal für die reine und angewandte Mathematik*, n. 94, 203-232.